\documentclass[11pt,a4paper,twoside]{article}

\usepackage{amsmath}
\usepackage{amssymb}
\usepackage{amsthm}
\usepackage[hang,flushmargin,bottom]{footmisc}
\usepackage[margin=0.35in,format=hang,font=small,labelfont=bf]{caption}
\usepackage{tikz}
\usepackage{chngcntr}
\usepackage[colorlinks=true,linkcolor=black,citecolor=black,urlcolor=black,filecolor=black]{hyperref}
\usepackage{palatino, mathpazo}

\usetikzlibrary{patterns}

\newcommand*\cleartoevenpage{
  \clearpage
  \ifodd\value{page}\hbox{}\thispagestyle{empty}\newpage\fi
}

\newcommand*\cleartooddpage{
  \clearpage
  \ifodd\value{page}\else\hbox{}\thispagestyle{empty}\newpage\fi
}

\newcounter{sect}

\newtheorem*{theorem*}{Theorem}

\newtheorem*{proposition*}{Proposition}

\counterwithin*{theorem}{sect}
\counterwithin*{equation}{sect}
\counterwithin*{figure}{sect}
\counterwithin*{table}{sect}

\newcommand{\talkTitle}[1]{\newpage\section*{\raggedright\textsf{\textmd{\LARGE #1}}}\hrule\stepcounter{sect}}
\newcommand{\talkPresenter}[1]{\textbf{\large #1}\vspace{12pt}}
\newcommand{\talkAffiliation}[1]{(#1)}

\newcommand{\subsect}[1]{\vspace{-6pt}\subsection*{#1}\vspace{-6pt}}
\newcommand{\subsectFirst}[1]{\vspace{-12pt}\subsection*{#1}\vspace{-6pt}}

\let\thebibliographyOld\thebibliography
\renewcommand\thebibliography[1]{
  \thebibliographyOld{#1}
  \setlength{\itemsep}{0pt plus 0.3ex}
}

\newenvironment{talkBibliog}[1][9]{\vspace{-6pt}\begin{thebibliography}{#1}\small}{\end{thebibliography}}

\begin{document}


\newenvironment{bullets} {\vspace{-9pt}\begin{itemize}\itemsep0pt} {\end{itemize}\vspace{-9pt}}

\newcommand{\+}{\hspace{0.07 em}}

\newcommand{\liminfty}[1][n]{\lim\limits_{#1\rightarrow\infty}}
\newcommand{\limsupinfty}[1][n]{\limsup\limits_{#1\rightarrow\infty}}
\newcommand{\liminfinfty}[1][n]{\liminf\limits_{#1\rightarrow\infty}}

\newcommand{\CCC}{\mathcal{C}}
\newcommand{\DDD}{\mathcal{D}}

\newcommand{\av}{\mathsf{Av}}
\newcommand{\pdash}{\text{--}}
\newcommand{\gr}{\mathrm{gr}}
\newcommand{\grup}{\overline{\gr}}
\newcommand{\grlow}{\underline{\gr}}
\newcommand{\red}{\mathrm{red}}
\newcommand{\inv}{\mathsf{inv}}
\newcommand{\Desc}{\mathsf{Des}}
\newcommand{\des}{\mathsf{des}}
\newcommand{\maj}{\mathsf{maj}}
\newcommand{\exce}{\mathsf{exc}}
\newcommand{\stat}{\mathsf{st}}

\newcommand{\plotpermnobox}[3][]  
{
  \foreach \y [count=\x] in {#3}
  {
    \ifnum0=\y {} \else {
      \fill[#1,radius=0.275] (\x,\y) circle;
    } \fi
  }
}

\newcommand{\plotpermgrid}[3][]  
{
  \foreach \x in {1,...,#2} \draw[very thin] (\x,.5)--(\x,#2.5) (.5,\x)--(#2.5,\x);
  \foreach \y [count=\x] in {#3} \fill[#1,radius=0.275] (\x,\y) circle;
  \draw[thick] (.5,.5) rectangle (#2.5,#2.5);
}

\newcommand{\plotpermmesh}[3][]  
{
  \foreach \x in {1,...,#2} \draw[thin] (\x,0.001)--(\x,#2.999) (0.001,\x)--(#2.999,\x);
  \foreach \y [count=\x] in {#3} \fill[#1,radius=0.2] (\x,\y) circle;
}

\talkTitle{Permutation patterns: basic definitions and notation}\label{sectGentleIntro}
\talkPresenter{David Bevan}
\talkAffiliation{The Open University}


\subsectFirst{Permutations, containment and avoidance}
A permutation is considered to be simply an arrangement of the numbers $1,2,\ldots,n$ for some positive $n$.
The length of permutation $\sigma$ is denoted $|\sigma|$, and
$S_n$ or $\mathfrak{S}_n$ is used for the set of
all permutations of length $n$. 

It is common to consider permutations graphically.
Given a permutation $\sigma=\sigma(1)\ldots\sigma(n)$, its \emph{plot} consists of the the points $(i,\sigma(i))$ in the Euclidean plane, for $i=1,\ldots,n$.

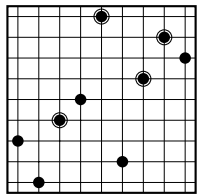
\begin{figure}[ht]
  $$
  \begin{tikzpicture}[scale=0.275]
    \plotpermgrid{9}{3,1,4,5,9,2,6,8,7}
    \draw [thin] (3,4) circle [radius=0.36];
    \draw [thin] (5,9) circle [radius=0.36];
    \draw [thin] (7,6) circle [radius=0.36];
    \draw [thin] (8,8) circle [radius=0.36];
  \end{tikzpicture}
  \vspace{-3pt}
  $$
  \caption{The plot of permutation $\mathbf{314592687}$ with a $\mathbf{1423}$ subpermutation marked}\label{figPermutation}
\end{figure}
A permutation, or \emph{pattern}, $\pi$ is said to be \emph{contained} in, or to be a \emph{subpermutation} of, another permutation $\sigma$, written $\pi\leqslant\sigma$ or $\pi\preccurlyeq\sigma$, if~$\sigma$ has a (not necessarily contiguous) subsequence whose terms are order isomorphic to (i.e.~have the same relative ordering as)
$\pi$.
From the graphical perspective,
$\sigma$ contains $\pi$ if the plot of $\pi$ results from erasing zero or more points from the plot of $\sigma$ and then rescaling the axes appropriately.
For example, $\mathbf{314592687}$ contains $\mathbf{1423}$
because the subsequence $\mathbf{4968}$ (among others) is ordered
in the same way as $\mathbf{1423}$ (see Figure~\ref{figPermutation}).

If $\sigma$ does not contain $\pi$, we say that $\sigma$ \emph{avoids} $\pi$.
For example, $\mathbf{314592687}$ avoids $\mathbf{3241}$ since it has no subsequence ordered in the same way as $\mathbf{3241}$.

If $\lambda$ is a list of distinct integers, the \emph{reduction} or \emph{reduced form} of $\lambda$, denoted
$\red(\lambda)$, is the permutation obtained from $\lambda$ by replacing
its $i$-th smallest entry with $i$. For example, we have $\red(\mathbf{4968}) = \mathbf{1423}$.
Thus, $\pi\leqslant\sigma$ if there is a subsequence $\lambda$ of $\sigma$ such that $\red(\lambda)=\pi$.

\subsect{Permutation structure}

\begin{figure}[ht]
  $$
  \begin{tikzpicture}[scale=0.275]
    \plotpermnobox{}{2,4,1,3,8,6,7,5}
    \draw[] (0.5,0.5) rectangle (4.5,4.5);
    \draw[] (4.5,4.5) rectangle (8.5,8.5);
  \end{tikzpicture}
  \qquad\qquad
  \begin{tikzpicture}[scale=0.275]
    \plotpermnobox{}{6,8,5,7,4,2,3,1}
    \draw[] (4.5,0.5) rectangle (8.5,4.5);
    \draw[] (0.5,4.5) rectangle (4.5,8.5);
  \end{tikzpicture}
  \qquad\qquad
  \begin{tikzpicture}[scale=0.275]
    \plotpermnobox{}{2,1,3,6,5,4,8,7}
    \draw[] (0.5,0.5) rectangle (2.5,2.5);
    \draw[] (2.5,2.5) rectangle (3.5,3.5);
    \draw[] (3.5,3.5) rectangle (6.5,6.5);
    \draw[] (6.5,6.5) rectangle (8.5,8.5);
  \end{tikzpicture}
  \vspace{-3pt}
  $$
  \caption{The direct sum $\mathbf{2413}\oplus\mathbf{4231}$, the skew sum $\mathbf{2413}\ominus\mathbf{4231}$, and the layered permutation $\mathbf{21}\oplus\mathbf{1}\oplus\mathbf{321}\oplus\mathbf{21}$}
  \label{figPermSums}
\end{figure}
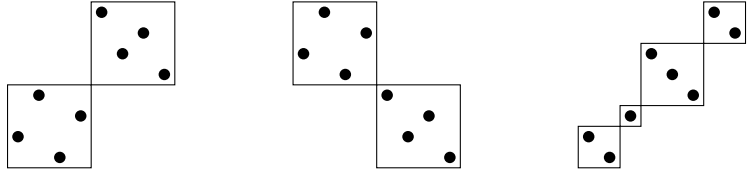

Given two permutations $\sigma$ and $\tau$ with lengths $k$ and $\ell$ respectively, their \emph{direct sum}\label{defDirectSum} $\sigma\oplus\tau$ is the permutation of length $k+\ell$ consisting of $\sigma$ followed by a shifted copy of $\tau$:
$$
(\sigma\oplus\tau)(i) \;=\;
\begin{cases}
  \sigma(i)   & \text{if~} i\leqslant k , \\
  k+\tau(i-k) & \text{if~} k+1 \leqslant i\leqslant k+\ell .
\end{cases}
$$
The \emph{skew sum} $\sigma\ominus\tau$ is defined analogously.
See Figure~\ref{figPermSums} for an illustration.

A permutation is called \emph{sum indecomposable}
if it cannot be expressed as the direct sum of two shorter permutations.
A permutation is \emph{skew indecomposable} if it cannot be expressed as the skew sum of two shorter permutations.
Every permutation has a unique representation as the direct sum of a sequence of sum indecomposable permutations, and also as the skew sum of a sequence of skew indecomposable permutations.
If a permutation is the direct sum of a sequence of \emph{decreasing} permutations,
then we say that the permutation is \emph{layered}.
See Figure~\ref{figPermSums} for an example.

An \emph{interval} of a permutation $\sigma$ corresponds to a contiguous sequence of indices $a,a+1,\ldots,b$ such
that the set of values $\{\sigma(i) : a\leqslant i \leqslant b\}$ is also contiguous.
Graphically,
an interval in a permutation is a square ``box'' that is not cut horizontally or vertically by any point
not in it.
Every permutation of
length $n$ has intervals of lengths $0$, $1$ and $n$. If a permutation $\sigma$ has no other intervals,
then $\sigma$ is said to be \emph{simple}.

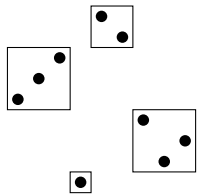
\begin{figure}[ht]
  $$
  \begin{tikzpicture}[scale=0.275]
    \plotpermnobox{}{5,6,7,1,9,8,4,2,3}
    \draw[] (0.5,4.5) rectangle (3.5,7.5);
    \draw[] (3.5,0.5) rectangle (4.5,1.5);
    \draw[] (4.5,7.5) rectangle (6.5,9.5);
    \draw[] (6.5,1.5) rectangle (9.5,4.5);
  \end{tikzpicture}
  \vspace{-3pt}
  $$
  \caption{The inflation $\mathbf{3142}[\mathbf{123}, \mathbf{1}, \mathbf{21}, \mathbf{312}] = \mathbf{567198423}$}
  \label{figPermInflat}
\end{figure}
Given a permutation $\sigma \in S_m$ and nonempty
permutations $\tau_1,\ldots,\tau_m$, the \emph{inflation} of $\sigma$ by $\tau_1,\ldots,\tau_m$, denoted
$\sigma[\tau_1,\ldots,\tau_m]$, is the permutation obtained by replacing each entry $\sigma(i)$
of $\sigma$ with an interval that is order isomorphic to $\tau_i$.
See Figure~\ref{figPermInflat} for an illustration.

A simple permutation is thus a permutation that cannot be expressed as
the inflation of a shorter permutation of length greater than $1$.
Conversely, every permutation except $\mathbf{1}$ is
the inflation of a unique simple permutation of length at least~$2$.

\begin{figure}[ht]
  $$
  \begin{tikzpicture}[scale=0.21]
    \plotpermnobox{14}{7,10,1,4,9,14,2,11,3,13,12,6,8,5}
    \draw[] (0.2,0.2) rectangle (14.8,14.8);
    \draw [thick] (1,7) circle [radius=0.4];
    \draw [thick] (2,10) circle [radius=0.4];
    \draw [thick] (6,14) circle [radius=0.4];
    \draw [thick] (3,1) circle [radius=0.4];
    \draw [thick] (7,2) circle [radius=0.4];
    \draw [thick] (9,3) circle [radius=0.4];
    \draw [thick] (14,5) circle [radius=0.4];
  \end{tikzpicture}
  \qquad\qquad\qquad
  \begin{tikzpicture}[scale=0.21]
    \plotpermnobox{14}{7,10,1,4,9,14,2,11,3,13,12,6,8,5}
    \draw[] (0.2,0.2) rectangle (14.8,14.8);
    \draw [thick] (1,7) circle [radius=0.4];
    \draw [thick] (3,1) circle [radius=0.4];
    \draw [thick] (6,14) circle [radius=0.4];
    \draw [thick] (10,13) circle [radius=0.4];
    \draw [thick] (11,12) circle [radius=0.4];
    \draw [thick] (13,8) circle [radius=0.4];
    \draw [thick] (14,5) circle [radius=0.4];
  \end{tikzpicture}
  \vspace{-3pt}
  $$
  \caption{The three left-to-right maxima and four right-to-left minima, and the two left-to-right minima and five right-to-left maxima, of a permutation}
  \label{figPermExtrema}
\end{figure}
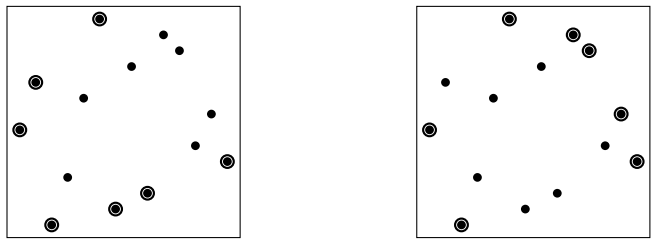
Sometimes we want to refer to the extremal points in a permutation.
A value in a permutation is called a \emph{left-to-right maximum} if it is larger than all the values to its left.
\emph{Left-to-right minima}, \emph{right-to-left maxima} and \emph{right-to-left minima} are defined analogously.
See Figure~\ref{figPermExtrema} for an illustration.

\subsect{Permutation statistics}

An \emph{ascent} in a permutation $\sigma$ is a position $i$ such that $\sigma(i)<\sigma(i+1)$.
Similarly, a \emph{descent} is a position $i$ such that $\sigma(i)>\sigma(i+1)$.
A pair of terms in a permutation $\sigma$ such that $i<j$ and $\sigma(i)>\sigma(j)$ is called an \emph{inversion}.

A permutation \emph{statistic} is simply a map from the set of permutations to the non-negative integers. Classical statistics include the following:
\begin{bullets}
  \item the number of descents \\ $\des(\sigma) = | \{i:\sigma(i)>\sigma(i+1)\} |$
  \item the number of inversions \\ $\inv(\sigma) = | \{(i,j):i<j \,\text{~and~}\, \sigma(i)>\sigma(j)\} |$
  \item the number of \emph{excedances} \\ $\exce(\sigma)=|\{i:\sigma(i)>i\}|$
  \item the \emph{major index}\footnote{\label{noteMacMahon}Named after Major Percy Alexander MacMahon.},
  the sum of the positions of the descents \\ $\maj(\sigma)=\sum_{\sigma(i)>\sigma(i+1)}i$
\end{bullets}

The statistics $\des$ and $\exce$ are equidistributed. That is, for all $n$ and $k$, the number of permutations of length $n$ with $k$ descents is the same as the number of permutations of length $n$ with $k$ excedances. Furthermore, $\inv$ and $\maj$ also have the same distribution.
Any permutation statistic that is distributed like $\des$ is said to be
\emph{Eulerian}, and a statistic that is distributed like $\inv$ is said
to be \emph{Mahonian}\footnote{See footnote~\ref{noteMacMahon}.}.

\subsect{Classical permutation classes}
The subpermutation relation is a partial order on the set of all permutations.
A classical \emph{permutation class}, sometimes called a \emph{pattern class}, is a set of permutations closed downwards (a {down-set}) under this partial order.
Thus, if $\sigma$ is a member of a permutation class $\CCC$ and $\tau$ is contained in $\sigma$, then it must be the case that $\tau$ is also a member of $\CCC$.
From a graphical perspective, this means that erasing points from
the plot of
a permutation in $\CCC$ always results in
the plot of
another permutation in $\CCC$
when the axes are rescaled appropriately.
It is common in the study of
classical permutation classes to
reserve the word ``class'' for
sets of permutations closed under taking subpermutations.

It is natural to define a classical permutation class ``negatively'' by stating the minimal set of permutations that it avoids.
This minimal forbidden set of patterns is known as the \emph{basis} of the class.
The class with basis $B$ is denoted $\av(B)$, and $\av_n(B)$ or $S_n(B)$ is used for the set of permutations of length $n$ in $\av(B)$.
As a trivial example, $\av(\mathbf{21})$ is the class of increasing permutations (i.e.~the identity permutation of each length).
As another simple example, the class of $\mathbf{123}$-avoiders, $\av(\mathbf{123})$, consists of those permutations that can be partitioned into two decreasing subsequences.

The basis of a permutation class is an {antichain} (a set of pairwise incomparable elements)
under the containment order,
and may be infinite.
Classes for which the basis is finite are called \emph{finitely based},
and those whose basis
consists of a single permutation are called \emph{principal} classes.

\subsect{Non-classical patterns}

Permutation patterns have been generalised in a variety of ways.

A \emph{barred} pattern is specified by a permutation with some entries barred ($\mathbf{5\bar{3}2\bar{1}4}$, for example).
If $\hat{\pi}$ is a barred pattern, let $\pi$ be the permutation obtained by removing all the bars in $\hat{\pi}$ ($\mathbf{5{3}2{1}4}$ in the example), and let $\pi'$ be the
permutation that is order isomorphic to the non-barred entries in $\hat{\pi}$ ($\mathbf{312}$ in the example).
An occurrence of barred pattern $\hat{\pi}$ in a permutation $\sigma$ is then an occurrence of $\pi'$ in $\sigma$ that is not part of an occurrence of $\pi$ in $\sigma$.
Conversely, for $\sigma$ to avoid $\hat{\pi}$, every occurrence in $\sigma$ of $\pi'$ must feature
as part of an occurrence of $\pi$.

\newcommand{\vinc}[1]{\mspace{1.2mu}\underline{\mspace{-1.2mu}#1\mspace{-1.2mu}}\mspace{1.2mu}}

A \emph{vincular} or \emph{generalised} pattern specifies adjacency conditions.
Two different notations are used. Traditionally, a vincular pattern is written as a permutation with dashes inserted between terms that
need not be adjacent and no dashes between terms that must be adjacent.
Alternatively, and perhaps preferably, terms that must be adjacent are underlined.
For example, $\mathbf{314265}$ contains two occurrences of $\mathbf{2\vinc{31}4}$ (or $\mathbf{2\pdash31\pdash4}$) and a single occurrence of $\mathbf{2\vinc{314}}$ ($\mathbf{2\pdash314}$),
but avoids $\mathbf{\vinc{23}\vinc{14}}$ ($\mathbf{23\pdash14}$).

A vincular pattern in which \emph{all} the terms must occur contiguously is known as a \emph{consecutive
pattern}.

In a \emph{bivincular} pattern, conditions are also placed on which terms must take adjacent \emph{values}.

\begin{figure}[ht]
  $$
  \begin{tikzpicture}[scale=0.4]
    \path [fill=gray!40] (0,2) rectangle (1,3);
    \path [fill=gray!40] (1,3) rectangle (2,5);
    \path [fill=gray!40] (4,2) rectangle (5,4);
    \plotpermmesh{4}{3,2,4,1}
  \end{tikzpicture}
  \vspace{-3pt}
  $$
  \caption{Mesh pattern $(\mathbf{3241}, \{(0, 2), (1, 3), (1, 4), (4, 2), (4, 3)\})$}\label{figMesh}
\end{figure}
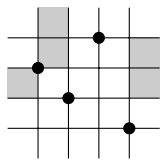

Classical, vincular and bivincular patterns are all example of the more general family of \emph{mesh} patterns.
Formally, a mesh pattern of length $k$ is a pair $(\pi,R)$ with $\pi\in S_k$ and $R\subseteq [0,k]\times[0,k]$ a set of pairs of integers.
The elements of $R$ identify the lower left corners of unit squares in the plot of $\pi$, which specify forbidden regions.
Mesh pattern $(\pi,R)$ is depicted by a figure consisting of the plot of $\pi$ with the forbidden regions shaded.
See Figure~\ref{figMesh} for an example.

An occurrence of mesh pattern $(\pi,R)$ in a permutation $\sigma$ consists of an occurrence of the classical pattern $\pi$ in $\sigma$ such that no elements of $\sigma$ occur in the shaded regions of the figure.
A vincular pattern is thus a mesh pattern in which complete columns shaded.

Sets of permutations defined by avoiding barred, vincular, bivincular or mesh patterns that are not closed under taking subpermutations are known as
\emph{non-classical} permutation classes.

\subsect{Growth rates}

Given a
permutation class $\CCC$,
we use $\CCC_n$ to denote the permutations of length $n$ in $\CCC$.
It is natural to ask how quickly the sequence
$(|\CCC_n|)_{n=1}^\infty$ grows.

In proving the Stanley--Wilf Conjecture, Marcos and Tardos
established that the growth of every
classical permutation class except the class of all permutations is at most exponential.
Hence, the \emph{upper growth rate} and \emph{lower growth rate} of a
class $\CCC$ are defined to be
$$
\grup(\CCC)=\limsupinfty{|\CCC_n|}^{1/n}
\qquad\text{~and~}\qquad
\grlow(\CCC)=\liminfinfty\+{|\CCC_n|}^{1/n} .
$$
The theorem of Marcos and Tardos states that $\grup(\CCC)$ and
$\grlow(\CCC)$ are both finite.

When $\grup(\CCC)=\grlow(\CCC)$, this quantity is called the \emph{proper growth
rate} (or just the \emph{growth rate}) of $\CCC$ and denoted $\gr(\CCC)$.
Principal classes, those of the form
$\av(\pi)$, are known to have proper growth rates.
The growth rate of $\av(\pi)$ is sometimes known as the \emph{Stanley--Wilf limit} of $\pi$ and denoted $L(\pi)$.
It is widely believed, though not yet proven, that every classical permutation class
has a proper growth rate.

\subsect{Wilf equivalence}

Given two classes, $\CCC$ and $\DDD$, one natural question is to determine whether they are equinumerous, i.e.~$|\CCC_n|=|\DDD_n|$ for every~$n$.
Two permutation classes that are equinumerous are said to be \emph{Wilf equivalent}
and the equivalence classes are called \emph{Wilf classes}.
If principal classes $\av({\sigma})$ and $\av({\tau})$ are Wilf equivalent, we simply say that $\sigma$ and $\tau$ are Wilf equivalent.

From the graphical perspective, it is clear that classes related by symmetries of the square are Wilf equivalent.
Thus, for example, $\av(\mathbf{132})$, $\av(\mathbf{231})$, $\av(\mathbf{213})$ and $\av(\mathbf{312})$ are equinumerous.
However, not all Wilf equivalences are a result of these symmetries. Indeed, as is well known, both $\av(\mathbf{123})$ and $\av(\mathbf{132})$ are counted by the Catalan numbers, so all permutations of length three are
in the same Wilf class.

\subsect{Generating functions}

The ordinary \emph{generating function} of a permutation class $\CCC$ is defined to be the formal power series
$$
C(z)
\;=\;
\sum_{n\geqslant0} |\CCC_n| z^n
\;=\;
\sum_{\sigma\in\CCC} z^{|\sigma|} .
$$
Thus, each permutation $\sigma\in\CCC$ makes a contribution of $z^{|\sigma|}$, the result being that, for each~$n$, the coefficient of $z^n$ is the number of permutations of length $n$.
Clearly, two classes are Wilf-equivalent if their generating functions are identical.

A generating function is \emph{rational} if it is the ratio of two polynomials.
A generating function $F(z)$ is \emph{algebraic} if it can be defined as the root of a polynomial equation.
That is, there exists a
bivariate polynomial $P(z, y)$ such that $P(z,F(z)) = 0$.

\begin{talkBibliog}

\bibitem{Bona2012}
Mikl{\'o}s B{\'o}na.
\newblock {\em Combinatorics of Permutations}.
\newblock Discrete Mathematics and its Applications. CRC Press, second edition,
  2012.

\bibitem{Kitaev2011}
Sergey Kitaev.
\newblock {\em Patterns in Permutations and Words}.
\newblock Springer, 2011.

\end{talkBibliog}

\end{document}